\documentclass[a4 paper]{amsart}

 \usepackage[english]{babel}
 \usepackage{latexsym}
 
 \newtheorem{theorem}{Theorem}[section]         
   \newtheorem{lemma}[theorem]{Lemma} 
              
 \newtheorem{corollary}[theorem]{Corollary}       
 \newtheorem{proposition}[theorem]{Proposition}

 \newtheorem{definition}[theorem]{Definition}     
                
 \newtheorem{example}[theorem]{Example} 
 \numberwithin{equation}{section}
 
 \newenvironment{namelist}[1]{%
 \begin{list}{}
     {
       
       \settowidth{\labelwidth}{#1}
       \setlength{\leftmargin}{1.1\labelwidth}
     }
   }{
 \end{list}}
 \newcommand{\nc}{\newcommand}
 \nc{\R}{\ensuremath{\mathbb{R}}}
 \nc{\sg}{\operatorname{sg}}
 \nc{\oC}{\boldsymbol{\ensuremath{\mathfrak C}}}
 \nc{\C}{\ensuremath{\mathfrak C}}
 \nc{\F}{\ensuremath{\mathbb{F}}}
 \nc{\M}{\ensuremath{\mathcal M}}
 \nc{\bM}{\ensuremath{\mathbf{M}}}
 \nc{\oM}{\boldsymbol{\M}}
 \title[mixed graphs and matroids]
 {a note on mixed graphs and matroids}
 \thanks{ 
 \emph{Secondary}: 05C20.
  }
 \author{J. Orestes Cerdeira and Raul Cordovil}
 \address{{}\newline
 Instituto Superior de Agronomia,\newline 
 Tapada da Ajuda,\newline 
 1349-017 Lisboa, Portugal}
 \email{orestes@isa.utl.pt}
 \address{{}\newline
 Departamento de Matem\'atica,\newline 
 Instituto Superior T\' ecnico \newline
  Av.~Rovisco Pais
  - 1049-001 Lisboa  - Portugal}
 \email{cordovil@math.ist.utl.pt}
 \thanks{The first author's research 
 was financially supported by the Portuguese Foundation for
 Science and Technology (FCT).
  The  second  author's research was 
 supported in part by FCT/POCTI/FEDER (Portugal) and
 the project FCT/FEDER/SAPIENS/36563/00.}
 \begin{document}
 \begin{abstract}
 {A mixed graph is a graph with some directed edges and 
 some undirected edges. 
 We introduce the notion of mixed matroids as a generalization
 of mixed graphs. A mixed matroid can be viewed as an oriented 
 matroid in which the signs over a fixed subset of the ground set 
 have been forgotten. We extend to mixed matroids standard 
 definitions from oriented matroids, establish basic properties, 
 and study questions regarding the reorientations of the unsigned 
 elements. In particular we address in the context of mixed matroids
 the $P$-connectivity and $P$-orientability issues which have been 
 recently introduced for mixed graphs.
 }
 \end{abstract}
 \maketitle 
 \section{Introduction}
 \label{Introduction}

 A mixed graph is a graph with a set of directed edges and a 
 set of undirected edges $A$.
 Recently, Arkin and Hassin \cite{arkin} consider the following 
 question regarding  the orientation of the edges of $A$.
 Given a set of ordered pair of vertices $P=(s_j,t_j), j=1,...,m$, 
 is there an orientation of the edges in $A$ such that in the resulting 
 digraph there is a (directed) $s_j\rightarrow t_j$ 
 path, for every $j=1,...,m$? 
 When such a $P$-\textit{orientation} exists the mixed graph is called
 $P$-\textit{orientable}. 
 
 Arkin and Hassin prove that recognizing whether a mixed graph has a 
 $P$-orientation is NP-\textit{complete}, and give a characterization
 of $P$-orientable graphs, for $|P|=2$, in terms of what they call
 $P$-essential edges. The mixed graph is $P$-\textit{connected} if, 
 for each pair $(s_j,t_j)\in P$, there is an orientation of the edges of 
 $A$ such that in the resulting digraph there is a (directed) 
 $s_j\rightarrow t_j$ path. An edge $e\in A$ is $P$-\textit{essential} if 
 the graph obtained from any of its two orientations is not $P$-connected.
 They prove the following.
 \begin{theorem} [Arkin and Hassin \cite{arkin}]\label{arkin_theorem}
 A $P$-connected mixed graph, with $|P|=2$, has a 
 $P$-orien\-ta\-tion if and only if it has no $P$-essential edges.
 \end{theorem}
 
 In this note we introduce the notion of mixed matroid as a 
 generalization of a mixed graph. 
 The concept of mixed 
 matroid, basic definitions and properties are given in section 
 \ref{mixed_matroids}. 
 In section \ref{P-orientation} we address the $P$-connectivity and 
 $P$-orientability questions for mixed matroids, and extend Theorem
 \ref{arkin_theorem}. 
 We end in section \ref{final_remark} with a remark
 concerning the possibility to obtain strong connected 
 digraphs from orientations of the undirected edges of 
 mixed graphs. 
 
 As an overall conclusion we trust that mixed matroids, and in 
 particular duality, give a consistent framework for the study of 
 issues such as connectivity in mixed graphs.

 \section{Mixed matroids}
 \label{mixed_matroids}
 
 We refer readers to \cite{Schrijver} as a 
 standard sources for graphs and matroids, 
 and \cite{book} for oriented matroids.
  
 A signed set (or vector) 
 $X,\, X\in [n]^{\{\pm1,0\}}$, is a (signed) mapping
 $$\sg_{X}: [n]\to \{\pm1,0\}.$$
 Define $X^+:=\{\ell\in [n]:\, \sg_{X}(\ell)= +1\}$, 
 $X^-:=\{\ell\in [n]:\, \sg_{X}(\ell)= -1\}$, 
 and let $\underline{X}:= X^+\cup X^-$ be the 
 support of the signed set $X$.
 We use the notation $\ell \in X$ meaning 
 that $\ell$ is in the support of $X$. 
 If $X$ and $Y$ are two signed sets, in order 
 to simplify notation,  
 $\underline{X}\cap \underline{Y}$ will be 
 denoted by $X\cap Y$.
 If $X$ is a signed set and $Y$ is a subset of $[n]$, 
 define $X\cap Y:=\underline{X}\cap Y$, and let
 $\sg_{X\setminus Y}$ be the map from $[n]$ 
 to $\{\pm1, 0\}$
  where  
 $$\sg_{X\setminus Y}(\ell)=\begin{cases}
 \sg_{X}(\ell)& \mbox{if}  
 \hspace{5mm} \ell\in X\setminus Y,\\
 0& \mbox{otherwise.} 
 \end{cases}
 $$
 Let $\M([n])$ [resp. $\oM([n])]$ denote a matroid
 [resp. oriented matroid] on the ground set $[n]$.
 To simplify
 we do not explicit refer the ground set
 $[n]$ simply writing $\M$ or $\oM$.  
 As usual $\oM^{\star}$ 
  denotes the  dual of $\oM$.  Finally, we denote by 
 $\C(\M)$ [resp. $\oC(\oM)$]  the set of circuits 
 [resp. signed circuits] 
 of $\M$ [resp. $\oM$]. 
 \begin{definition}
 {\em
 A \emph{mixed matroid} is a pair $\bM=(\oM, A)$, where 
 $\oM$ is an oriented matroid,  $A$ is a subset
 of the ground set $[n]$, and to which we associate the set
 $$\oC(\bM):=\{\sg_{C\setminus A}:\,\, \sg_{C}\in \oC(\oM)\},$$  
 that will be called the set
 of the \emph{circuits} $\bM$. 
 }
 \end{definition}
 It is useful to interpret $\bM$ as an oriented matroid whose signs 
 of a fixed subset $A$ have been forgotten in every signed circuit. 
 The set $A$ will be referred as the set of 
 \emph{unsigned} elements of the mixed matroid $\bM$.  
 Note that in the definition above we fixed 
 an orientation $\oM$ of the underlying
 (non-oriented) matroid because,
 as opposed to what happen with graphs,
 an orientable matroid may have different 
 orientations not in the same class of 
 reorientations.

 Examples of mixed matroids are obtained in the extreme cases 
 $A=\emptyset$ and $A=[n]$ where $\bM$ is, respectively, an 
 oriented matroid and an orientable matroid.
 
 \begin{definition}
 {\em
 Any oriented matroid $\oM^\prime$ on the class of reorientations 
 of $\oM$ over $A$ is said to be \emph{coherent} with the 
 mixed matroid $\bM=(\oM, A)$. 
 }
 \end{definition}
 Consider two mixed matroids $\bM=(\oM, A)$ and 
 $\bM^{\prime}=(\oM^{\prime}, A^{\prime})$ with the same underlying matroid.
 If $A^{\prime} \subset A$, and $\oM^{\prime}$ is coherent with $\bM$,
 we say that $\bM^{\prime}$ is obtained from $\bM$ by a 
 permissible signature of the elements of $A\setminus A^\prime$, and write
 $\bM^{\prime}:=\bM(\overrightarrow{A\setminus A^\prime})$. 

 The operations  \emph{deletion} and \emph{contraction}  on oriented 
 matroids naturally extend to
 a mixed oriented matroid $\bM=(\oM, A)$. 
 The \emph{deletion} $\bM\setminus X$ is
 the mixed matroid  $(\oM\setminus X, A\setminus X)$ where 
 $$\oC(\bM\setminus X):=\{\sg_{C\setminus A}:\,\,  
 \sg_{C}\in \oC(\oM\setminus X)\}.$$  
 The \emph{contraction} $\bM/X$ is
 the mixed matroid  $(\oM/X, A\setminus X)$\, where
 $$\oC(\bM/X):=\{\sg_{C\setminus A}:\,\, \sg_{C}\in \oC(\oM/X)\}.$$  
 
The mixed matroids support a notion of duality 
 similar to that of oriented matroids.
 The \emph{dual} or \emph{orthogonal} of  $\bM$ 
 is the mixed matroid $\bM^{\star}=(\oM^{\star},A)$ 
 whose set of circuits is 
 $$
 \oC^{\star}(\bM):=
 \{ \sg_{C^{\star}\setminus A}:\,\,  
 \sg_{C^{\star}}\in \oC(\oM^{\star})\}.$$
 \begin{example}
 {\em 
 Consider a mixed graph $G$ with a set $A$ of undirected edges. 
 Fix labels $e_1, e_2,\dots, e_n$ on the edges of $G$, and 
 identify the subsets of 
 edges with the subsets of indices of the corresponding 
 labels.
 The graphic (non oriented) matroid $\M(G)$ is a matroid on 
 the ground set $[n]$ and
 such that $C$ is a circuit of $\M(G)$ if and only if it 
 is a circuit of the graph $G.$
 The partial signatures of the edges of $G$ determine 
 a mixed graphic matroid.
 Indeed every partially signed circuit
 $C$ [resp. signed cocircuit or signed minimal cut 
 set $C^{\star}$] of the graph $G$  
 fixes a pair of opposite maps  $\pm\sg_{C\setminus A}$
 [$\pm\sg_{C^{\star}\setminus A}$].  
 Thus, the mixed graph $G$ determines a pair of dual 
 mixed oriented matroids
 $\bM(G)$ and $\bM^{\star}(G)$. 
 It is clear that every  possible orientation of the edges 
 of $A$ determines a digraph 
 $\widetilde{G}$, and the oriented matroid $\oM(\widetilde{ G})$
 and its dual $\oM^{\star}(\widetilde{G})$ are coherent with 
 the mixed matroids
 $\bM(G)$ and $\bM^{\star}(G)$, respectively.
 }
 \end{example}
 The following are dual propositions. 
 \begin{proposition}\label{pro} Let $\bM=(\oM, A)$ be 
 a coloop free mixed  matroid. 
 The two conditions are equivalent:
 \begin{namelist}{xxxx}
 \item[~~$\circ$] $\bM/A$ is a totally cyclic oriented matroid; 
 \item[~~$\circ$] There is a totally cyclic reorientation 
 of \,$\oM$ coherent with $\bM$. 
 \end{namelist}
  \end{proposition}
 \begin{proposition}\label{pro2} Let $\bM=(\oM, A)$ be a 
 loop free mixed  matroid. 
 The two conditions are equivalent:
  \begin{namelist}{xxxx}
 \item[~~$\circ$] $\bM\setminus A$ is an acyclic oriented matroid. 
 \item[~~$\circ$] There is an acyclic  reorientation 
 of \,$\oM$ coherent with $\bM$. \qed
 \end{namelist}
  \end{proposition}
 \noindent{\em Proof of Proposition  \ref{pro}.} 
  If a reorientation of $\oM$  coherent with $\bM$ 
 is totally cyclic, 
 then it is a consequence of the definitions that
 $\bM/A$ is also  a totally cyclic oriented matroid.
  
 Suppose that $\bM/A$ is a totally cyclic oriented 
 matroid, and $\oM$ is not.
 Given that the contraction of any element of a coloop free 
 matroid does not creates coloops, 
 we can consider without loss of generality that  $A=\{e\}$.
 Since $\oM^{\star}$ is not acyclic and $\oM$ is  coloop free, 
 there is a positive cocircuit $D$ with at least two elements that  
 includes the element $e$. 
 By hypothesis $\oM^{\star}\setminus e$
 is an acyclic oriented matroid, and the elimination axiom 
 for cocircuits ensures that $\oM^{\star}$
 has no signed cocircuit $D_{1}$
 with $D_{1}^{-}=\{e\}$.
 Indeed, the existence of $D$ and $D_1$ together would imply the 
 existence of a positive cocircuit not containing $e$, 
 which is a contradiction.
 It then follows that  ${}_{-\{e\}}\oM^{\star}$ is acyclic, and 
 by duality $\widetilde{\oM}:={}_{-\{e\}}\oM$ is a 
 totally cyclic oriented matroid
 coherent with $\bM$.\qed\par

 \section{P-orientation}
 \label{P-orientation}
 
 This section extends Theorem \ref{arkin_theorem} to matroids. 
 
 Let $\bM=(\oM, A)$ be a mixed matroid. The circuit $C\in \oC(\bM) $ is 
 said to be \emph{positive} if 
 $\sg_{C\setminus A}(C\setminus A)=\{+1\}$. 
 
 The following two definitions are generalizations 
 (with a slight modification) 
 of the definitions introduced by Arkin and Hassin \cite{arkin}, and 
 which we presented in section 
 \ref{Introduction}. 
 
 Let $P$ be a fixed set of signed elements of the mixed matroid $\bM$.
 \begin{definition}\label{connectivity}
 {\em
 The mixed matroid $\bM$ is $P$-\emph{connected} if the two 
 following conditions are verified:
 \begin{namelist}{xxxxxxxxxx}
 \item[~~$(\ref{connectivity}.1)$] For every element $p\in P$, there is a
 positive circuit $C$ of $\bM$ such that \newline $C\cap P=\{p\}$;
 \item[~~$(\ref{connectivity}.2)$] $\bM$ is 
  totally cyclic.
 \end{namelist}  
 }
 \end{definition}
 The condition $(\ref{connectivity}.2)$ is technical and 
 not essential to the 
 definition. Without it the main results would still hold, 
 however, the proofs would be more involved.

 Note that if $A$ is a basis of $\M$ and $P= [n]\setminus A$, 
 then $\bM=(\oM, A)$ is clearly $P$-connected. 
 
 \begin{definition}\label{orientable}
 {\em
 The mixed matroid $\bM$ has a $P$-\emph{orientation} if there 
 exists a $P$-connected oriented matroid coherent with $\bM$.
 }
  \end{definition}
 \begin{lemma}
 Let $\bM=(\oM, A)$ be a $P$-connected mixed matroid. Then $P$ is 
 an independent set of the matroid $\oM^\star$. 
 Moreover, for every $a\in A$, $\bM/a$ is also a
 $P$-connected matroid. If in addition $\bM$ is $P$-orientable,
 then $\bM/a$ is also $P$-orientable.
 \end{lemma}
 \begin{proof} 
 Note that if $\bM$ is $P$-connected the condition  
 $(\ref{connectivity}.1)$ implies, by the
 orthogonality property of circuits and cocircuits, that $P$ 
 is an independent set of the matroid $\oM^\star$. 
 The other results are direct consequences of the definitions.
 \end{proof}
 
 We can give to Definition~\ref{orientable} a more 
 geometrical interpretation.
 Recall that a set (hyperplane) $F\subset [n]$ is 
 a \emph{facet} of the oriented matroid $\oM$ 
 if there is a positive
 cocircuit  $C\in \oC(\oM^\star)$ 
 such that $F=[n]\setminus C$. 
 \begin{proposition}\label{str}
 Let $\bM$ be a totally cyclic mixed $P$-connected matroid, and 
 suppose that $P$ has at  least two elements. The two conditions 
 are equivalent: 
 \begin{namelist}{xxxxxxxxxx}
 \item[~~$(\ref{str}.1)$] $\bM$  has a $P$-orientation;
 \item[~~$(\ref{str}.2)$] There is a totally cyclic oriented 
 matroid $\oM$ coherent with $\bM$ such that,
 for every $p\in P$, the set  $P\setminus p$ is on a facet of 
 $\oM^\star$.
 \end{namelist}
 \end{proposition}
 \begin{proof}
 This is the dual of Definition~\ref{orientable}.
 \end{proof}
 
 The definition below extends the definition in   
 \cite{arkin} of $P$-essential edges of a mixed graph.
 \begin{definition}\label{excluded}
 {\em
 Let $\bM$ be a $P$-connected mixed matroid, and $e$ an unsigned
 element of $\bM$. The element $e$ is $P$-\emph{essential} 
 if the two possible opposite permissible signatures 
 $\bM(\vec{e})$ and ${}_{-\{e\}}\bM(\vec{e})$ are 
 not $P$-connected.
 }
 \end{definition}
 
 The following technical proposition is fundamental in this section. 
 See the mixed graph $G$ of Figure 1 for an illustration of the 
 result for mixed graphs. 
 \begin{proposition}\label{equiv} Let $\bM$ be a $P$-connected mixed 
 matroid. Suppose that $P=\{p_1, p_2\}$ and neither $p_1$ nor $p_2$ 
 is a loop. The conditions $(\ref{equiv}.1)$ and $(\ref{equiv}.2)$ are
 equivalent:
 \begin{namelist}{xxxxxxxx}
 \item[~~$(\ref{equiv}.1)$] The unsigned element $e$ is 
 $P$-\emph{essential};
 \item[~~$(\ref{equiv}.2)$]
 \begin{namelist}{xxxx}
 \item[~~$a)$] 
 If $C$ is a positive circuit and $|C\cap P| =1$ then $e\in C$.
 \item[~~$b)$] If $C_1$ and $C_2$ are two positive circuits with 
 $p_i \in C_i,\, i=1,2$, there is positive circuit $C_{12}$ 
 containing $p_1$ and $p_2$ such that 
 $C_{12}\subset C_1\cup C_2\setminus e$.
 \end{namelist}
 \end{namelist}
 \end{proposition}
 \begin{proof}
 We first prove that  $(\ref{equiv}.1)\Longrightarrow 
 (\ref{equiv}.2)$. 
 The condition $(\ref{equiv}.2a)$ is a direct 
 consequence of the cardinality 
 of $P$. 
 To prove $(\ref{equiv}.2b)$ consider two  
 positive circuits $C_1$ and $C_2$ of $\bM$,
 such that $C_i\cap P= \{p_i\},\, i=1,2$.
 Let $\bM(\vec{e})$ be the mixed matroid 
 obtained from $\bM$ by a permissive signature 
 of the element $e$, and 
 let $C_i(\vec{e}),\, i=1,2,$ be the circuits of 
 $\bM(\vec{e})$ corresponding to $C_i$. 
 As $e$ is $P$-essential it necessarily has 
 opposite signs in the circuits $C_1(\vec{e})$ and
 $C_2(\vec{e})$. The circuit elimination axiom
 ensures that there is a circuit $C_{12}$ 
 containing $p_1$ and 
 contained in $C_1\cup C_2\setminus e$. 
 The signed circuit axiom 
 ensures that  $C_{12}$ is a positive circuit 
 of the mixed matroid $\bM$. 
 Since $p_1 \in C_{12}$ and $e\not \in C_{12}$, we 
 conclude from 
 $(\ref{equiv}.2a)$ that necessarily 
 $p_2\in C_{12}$ and 
 that $(\ref{equiv}.2b)$ is true.
 
 We now prove that 
 $(\ref{equiv}.2)\Longrightarrow 
 (\ref{equiv}.1)$. 
 Suppose that $e$ is not a $P$-essential element of
 $\bM$, and let $C_{1}(\vec{e})$ and $C_{2}(\vec{e})$ 
 be positive circuits of a permissible 
 signature $\bM(\vec{e})$, such that 
 $C_{i}(\vec{e})\cap P=\{p_i\}$ and $e\in C_{i}(\vec{e})$, $i=1,2$.
 The condition $(\ref{equiv}.2b)$ ensures the 
 existence of 
 a positive circuit
 ${C}_{12}$ of $\bM(\vec{e})$ such that 
 $$p_{1}, p_{2} \in {C}_{12} 
 \subset {C_{1}(\vec{e})}\cup
 {C_{2}(\vec{e})}\setminus e.$$
 Note that as $e\not \in C_{12}$ this  is also a positive circuit of $\bM$.
 By the elimination of $e$ in  $C_{1}(\vec{e})$ and 
 $C_{2}(\vec{e})$ 
 it follows from  the signed elimination axiom that  $p_{1}$
 and $p_{2}$ have opposite signs in the 
 signed circuit $C_{12}$, which is a contradiction.
 Hence, $e$ is a $P$-essential element of $\bM$.
 \end{proof}
 \begin{figure}
 \begin{center}
 \thicklines
 \begin{picture}(-70,120)(100,10)
 \put(70,10) {\vector(-1,0){67}}
 \put(130,10) {\vector(-1,0){57}}
 \put(0,10) {\circle*{6}}\put(-15,40){$p_1$}
 \put(70,10) {\circle*{6}}
 \put(130,10) {\circle*{6}}\put(135,40){$p_2$}
 \put(0,10) {\vector(0,1){57}}
 \put(0,70) {\circle*{6}}
 \put(0,70){\vector(1,0){67}}
 \put(70,70){\vector(1,0){57}}
 \put(70,70) {\circle*{6}}\put(60,40){$e$}
 \put(70,10){\line(0,1){60}}
 \put(130,70) {\vector(0,-1){57}}
 \put(130,70) {\circle*{6}}
 \end{picture}
 \caption{The mixed graph $G$.}
 \end{center}
  \end{figure}
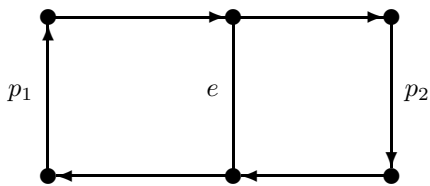
 
 We can now generalize  
 Theorem \ref{arkin_theorem}. 
 We highlight that the definition of $P$-essential elements 
 given by the Proposition~\ref{equiv}
 simplifies 
 the proof in \cite{arkin}.
 \begin{theorem}
 A $P$-connected mixed matroid $\bM=(\oM, A)$, with $|P|=2$, 
 has a $P$-orien\-tation if and only if it has no $P$-essential 
 elements.
 \end{theorem}
 \begin{proof} 
 The necessary condition is clear. To prove that it is 
 sufficient we use induction on the cardinal of $A$. 
 For $|A|= 0$ there is nothing to prove.  
 If $A=\{e\}$, the result follows directly from the hypothesis 
 that $e$ is not $P$-essential.
 Suppose now that the result is true for all sets of unsigned  
 elements with cardinality at most $m-1$, and let $2 \leq|A|=m$. 
 Since contracting any unsigned element of a mixed oriented 
 matroid with no $P$-essential elements does not create 
 $P$-essential elements, the induction hypothesis implies that, 
 for every $e\in A$, $\bM/e$ is $P$-orientable. 
 Hence, there is a signature of the elements of $A\setminus e$ 
 such that the mixed matroid 
 $\bM^\prime:=\bM(\overrightarrow{A\setminus e})$ is $P$-connected. 
 A contradiction results if $\bM^\prime$ is not 
 $P$-orientable, i.e., the element $e$ is $P$-essential in 
 $\bM^\prime$. In fact, let $P=\{p_1, p_2\}$.
 Proposition~\ref{equiv} ensures that, for every pair of positive 
 circuits $C_1$ and $C_2$ of $\bM^\prime$ with  
 $p_i \in C_i,\, i=1,2$, necessarily $e\in C_1\cap C_2$, and 
 there is a positive circuit $C_{12}$ such that 
 $p_1,p_2 \in C_{12} \subset C_1 \cup C_2 \setminus e$. 
 We can now forget about the signs over $A$, and use 
 Proposition~\ref{equiv} to conclude that element $e$ 
 is also $P$-essential in $\bM$, which is a contradiction. Thus,
 $\bM$ has a $P$-orientation.
 \end{proof}

 \section{final remark}
 \label{final_remark}
 
 We point out that the graph version of Proposition \ref{pro} 
 is a well known result 
  of mixed graphs.
 
 Observe that coloop free matroids correspond to graphs where each 
 component is 2-edge connected,
 and totally cyclic oriented matroids correspond to digraphs 
 where each component is strongly connected. 
 Thus, it follows from the 
 Proposition \ref{pro} that a 2-edge connected mixed graph
 has an orientation of the undirected edges so that 
 the resulting digraph is strongly connected if and only if 
 the digraph obtained contracting the undirected edges is strongly
 connected. 
 Since 2-edge connectivity is clearly a necessary condition to turn 
 a mixed graph into a strongly connected digraph, this is precisely 
 the following well known result. 
 \begin{theorem} [\cite{Schrijver}, Theorem 61.4]\label{main}
 Let $G=(V,E)$ be a graph in which part of the edges is oriented 
 (a mixed graph). Then the remainder of the edges can be oriented 
 so as to obtain a strongly connected digraph 
 if and only if $G$ is 2-edge connected and there is no non empty 
 proper subset $U$ of\, $V$ such that all edges
 in $\delta (U)$ are oriented from $U$ to $V\setminus U.$ \qed
 \end{theorem}
 \par
 From the Proposition \ref{pro2} we can establish the dual of 
 Theorem \ref{main}.
 \begin{corollary} Let $G=(V,E)$ be a 
  mixed graph. Then the non oriented  edges can be oriented 
 so as to obtain an acyclic  digraph 
 if and only  for every     oriented edge $e$ there is a non empty 
 proper subset $U$ of\, $V$ such that $e\in \delta (U)$  and all edges
 in $\delta (U)\setminus e$ are non oriented or oriented from $U$ to $V\setminus U$ 
 and $e\in \delta (U)$.
 \qed
 \end{corollary}
 \vspace{2ex}
 
 \centerline{\textbf {Acknowledgements}}
 \vspace{1ex}
 \noindent We are grateful to Paulo Barcia and David Forge  for assistance and discussion.
 \vspace{2ex}

 
 \end{document}